\documentclass[letterpaper, 10 pt, conference]{ieeeconf}  % Comment this line out
\IEEEoverridecommandlockouts        % This command is only
\title{\LARGE \bf Cyber-Physical Attacks in Power Networks:\\
  Models, Fundamental Limitations and Monitor Design}

\author{Fabio Pasqualetti, Florian D\"orfler, and Francesco Bullo%
  \thanks{Fabio Pasqualetti, Florian D\"orfler, and Francesco Bullo are
    with the Center for Control, Dynamical Systems and Computation,
    University of California at Santa Barbara, {\tt
      \{fabiopas,dorfler,bullo\}@engineering.ucsb.edu}}%
}

\usepackage{graphicx,color}
\usepackage{amsmath}
\usepackage{amssymb}
\usepackage{mathrsfs}
\usepackage[vlined,ruled]{algorithm2e}
\usepackage{subfigure}

\newtheorem{theorem}{Theorem}[section]
\newtheorem{lemma}[theorem]{Lemma}
\newtheorem{definition}{Definition}

\newtheorem{remark}{Remark}

\newcommand{\subscr}[2]{{#1}_{\textup{#2}}}
\newcommand{\until}[1]{\{1,\dots,#1\}}

\newcommand{\Ker}{\operatorname{Ker}}

\newcommand{\Image}{\operatorname{Im}}

\newcommand{\Star}{\mathcal{S}^*}

\newcommand{\real}{\mathbb{R}}
\newcommand{\complex}{\mathbb{C}}
\newcommand{\transpose}{\mathsf{T}} %or \top or \intercal

\newcommand{\Basis}{\operatorname{Basis}}

\newcommand{\mc}{\mathcal}

\newcommand\oprocendsymbol{\hbox{$\square$}}
\newcommand\oprocend{\relax\ifmmode\else\unskip\hfill\fi\oprocendsymbol}

\renewcommand{\theenumi}{(\roman{enumi})}

\newenvironment{pfof}[1]{\vspace{1ex}\noindent{\itshape Proof of
    #1:}\hspace{0.5em}} {\hfill\QED\vspace{1ex}}

\newenvironment{pf}{\vspace{1ex}\noindent{\itshape Proof:}\hspace{0.5em}}
{\hfill\QED\vspace{1ex}}

\renewcommand{\baselinestretch}{0.995}
\begin{document}
\maketitle
\thispagestyle{empty}
\pagestyle{empty}

%%%%%%%%%%%%%%%%%%%%%%%%%%%%%%%%%%%%%%%%%%%%%%%%%%%%%%%%%%%%%%%%%%%%%%%%%%%%%%%%

\begin{abstract}
  Future power networks will be characterized by safe and reliable
  functionality against physical malfunctions and cyber attacks. This
  paper proposes a unified framework and advanced monitoring
  procedures to detect and identify network components malfunction or
  measurements corruption caused by an omniscient adversary. We model
  a power system under cyber-physical attack as a linear
  time-invariant descriptor system with unknown inputs. Our attack
  model generalizes the prototypical stealth, (dynamic) false-data
  injection and replay attacks. We characterize the fundamental
  limitations of both static and dynamic procedures for attack
  detection and identification. Additionally, we design
  provably-correct (dynamic) detection and identification procedures
  based on tools from geometric control theory.  Finally, we
  illustrate the effectiveness of our method through a comparison with
  existing (static) detection algorithms, and through a numerical
  study.
\end{abstract}

%%%%%%%%%%%%%%%%%%%%%%%%%%%%%%%%%%%%%%%%%%%%%%%%%%%%%%%%%%%%%%%%%%%%%%%%%%%%%%%%

\section{Introduction}
\label{sec:introduction}
\textbf{Problem setup} Recent studies and real-world incidents have
demonstrated the inability of the power grid to ensure a reliable service
in the presence of network failures and possibly malignant actions
\cite{ARM-RLE:10,MA:11}.
%JM:07
Besides failures and attacks on the {\em physical} power grid
infrastructure, the envisioned future smart grid is also prone to {\em
  cyber} attacks on its communication layer. In short, cyber-physical
security is a fundamental obstacles challenging the smart grid vision.

A classical mathematical model to describe the grid on the
transmission level is the so-called {\em structure-preserving power
  network model}, which consists of the dynamic {\em swing equation}
for the generator rotor dynamics, and of the algebraic
\emph{load-flow} equation for the power flows through the network
buses \cite{PWS-MAP:98}. In this work, we consider the linearized
small signal version of the structure-preserving model, which is
composed by the linearized swing equation and the {\em DC power flow
  equation}.  The resulting linear continuous-time descriptor model
%\cite{PK-VLM:06}
of a power network has also been studied for estimation and security
purposes in \cite{ES:04,ADDG-ST:10,FP-AB-FB:10u}.

\begin{figure}
    \centering
    \includegraphics[width=.9\columnwidth]{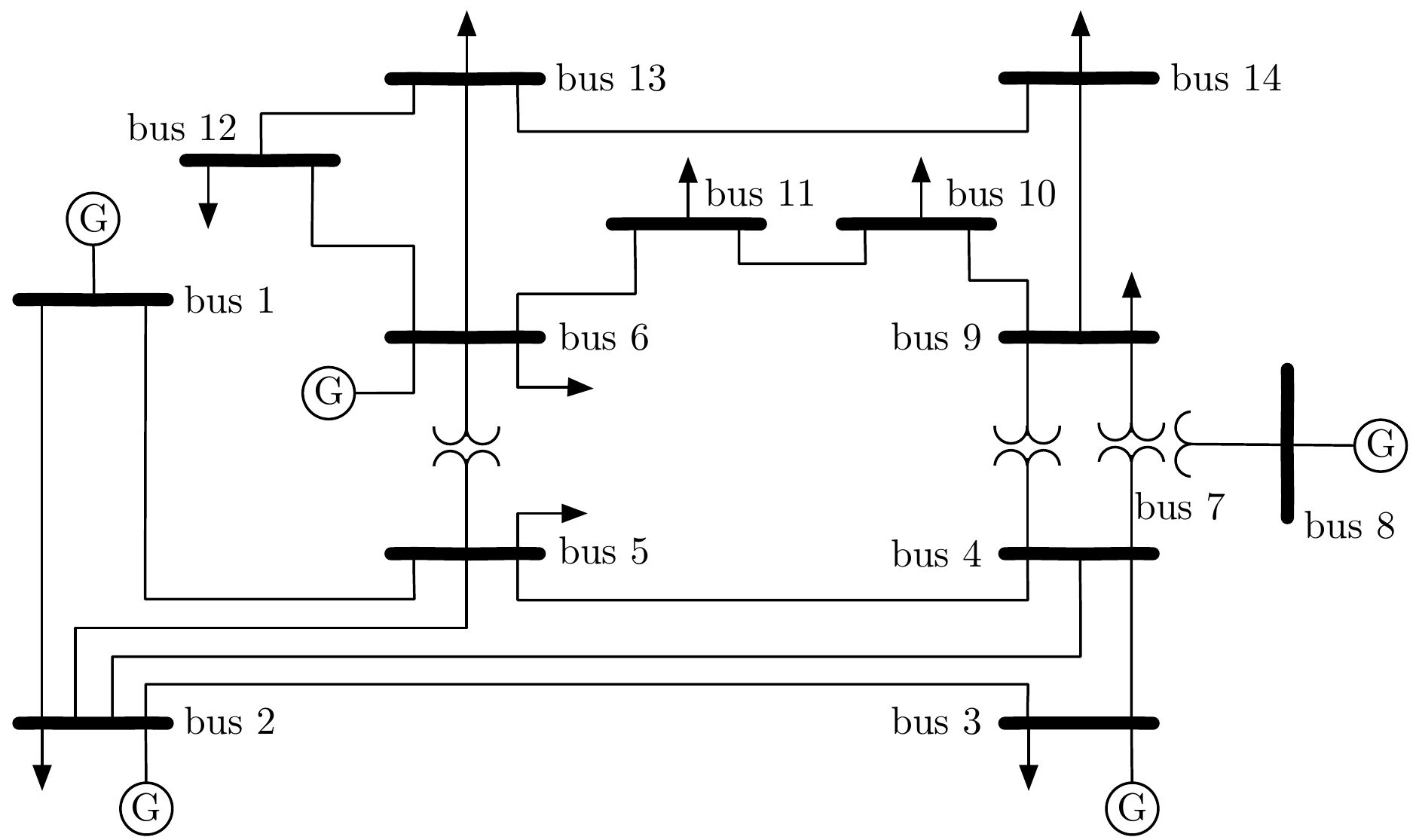}
    \caption{For the here represented IEEE 14 bus system, if the
      voltage angle of one bus is measured exactly, then a cyber
      attack against the measurements data is always detectable by our
      dynamic detection procedure. In contrary, as shown in
      \cite{YL-MKR-PN:09}, a cyber attack may remain undetected by a
      static procedure if it compromises as few as four measurements.}
    \label{fig:ieee14}
\end{figure}

\textbf{From static to dynamic detection} Existing approaches to
security and stability assessment are mainly based upon static
estimation techniques for the set of voltage angles and magnitudes at
all system buses, e.g., see \cite{AA-AGS:04}. Limitations of these
techniques have been often underlined, especially when the network
malfunction is intentionally caused by an omniscient attacker 
% aiming at affecting the network functionality while remaining
% undetected
\cite{YL-MKR-PN:09,AT-AS-HS-KHJ-SSS:10}. The development of security
procedures that exploit the dynamics of the power network is
recognized~\cite{NB-TB-AB-VB-GC-DC-AF-LF-MGL-BFW-JNW:92} as an
outstanding important problem. We remark that the use of static state
estimation and detection algorithms has been adopted for many years
for several practical and technological reasons. First, because of the
low bandwidth of communication channels from the measuring units to
the network control centers, continuous measurements were not
available at the control centers, so that the transient behavior of
the network could not be captured. Second, a sufficiently accurate
dynamic model of the network was difficult to obtain or tune, making
the analysis of the dynamics even harder. As of today, because of
recent advances in hardware technologies, e.g., the advent of {\em
  Phasor Measurement Units} and of large bandwidth communications, and
in identification techniques for power system parameters
\cite{AC-JHC-AS:11}, these two limitations can be overcome. Finally, a
dynamic estimation and detection problem was considered much harder
than the static counterpart. We address this theoretic limitation by
improving upon results presented in \cite{FP-AB-FB:09b,SS-CH:10a} for
the security assessment of discrete time dynamical networks.

\textbf{Literature review on dynamic detection} Dynamic security has
been approached via heuristics and expert systems, e.g., see
%SM-KEB:97,KW-AE-RA:97
\cite{KYL:08}.  Shortcomings of these methods include reliability and
accuracy against unforeseen system anomalies, and the absence of
analytical performance guarantees. A different approach relies on
matching a discrete-time state transition map to a series of past
measurements via Kalman filtering, e.g., see
%PR-TVC:90,GDP-SST:98
\cite{ACZDS-JCSDS-AMLS:02,UAK-MDI-JMFM:10} and the references
therein. Typically, these transition maps are based on heuristic
models fitted to a specific operating point
%PR-TVC:90,GDP-SST:98
\cite{ACZDS-JCSDS-AMLS:02}. Clearly, such a pseudo-model poorly
describes the complex power network dynamics and suffers from
shortcomings similar to those of expert systems methods.
In \cite{UAK-MDI-JMFM:10}, the state transition map is chosen more
accurately as the linearized and Euler-discretized power network
dynamics. The local observability of the resulting linear
discrete-time system is investigated in \cite{UAK-MDI-JMFM:10}, but in
the absence of unforeseen attacks.
Finally, in \cite{DK-XF-SL-TZ-KLBP:10} a graph-theoretic framework is
proposed to evaluate the impact of cyber attacks on a smart grid and
empirical results are given. 

Recent approaches to dynamic security consider continuous-time power
system models and apply dynamic techniques
\cite{ES:04,ADDG-ST:10,FP-AB-FB:10u,AT-HS-KHJ:10}. While
\cite{AT-HS-KHJ:10} adopts an overly simplified model neglecting the
algebraic load flow equations, the references
\cite{ES:04,ADDG-ST:10,FP-AB-FB:10u} use a more accurate network
descriptor model.
% We remark that
% \cite{AT-HS-KHJ:10} uses an overly simplified and even ill-posed model
% neglecting the algebraic load flow equations. On the contrary, the
% work in~\cite{ES:04,ADDG-ST:10,FP-AB-FB:10u} uses the accurate
% descriptor system model.
In \cite{ADDG-ST:10} different failure modes are modeled as instances
of a switched system and identified using techniques from hybrid
control. This approach, though elegant, results in a severe
combinatorial complexity in the modeling of all possible
attacks. % Our earlier work \cite{FP-AB-FB:10u}, instead, relies on the
% notion of input identifiability and proposes a structured system
% approach to security assessment. Specifically,
In our earlier work \cite{FP-AB-FB:10u}, under the assumption of
generic network parameters, we state necessary and sufficient
conditions for identifiability of attacks based on the network
topology. Finally, in \cite{ES:04} dynamical filters are designed to
isolate certain predefined failures of the network components. With
respect to this last work, we assume no a priori knowledge of the set
of compromised components and of their compromised behavior. Our
results generalize and include those of \cite{ES:04}.

\textbf{Contributions} This paper's contributions are fourfold. %
First, we provide a unified modeling framework for dynamic power networks
subject to cyber-physical attacks. For our model, we define the notions of
{\em detectability} and {\em identifiability} of an attack by its effect on
output measurements.  Informed by the classic work on geometric control
theory~\cite{WMW:85,HLT-AS-MH:01}, our framework includes the
\emph{deterministic static detection problem} considered in
% FCS-JW:70
\cite{YL-MKR-PN:09,AT-AS-HS-KHJ-SSS:10}, and the prototypical
\emph{stealth} \cite{DG-HS:10}, \emph{(dynamic) false-data injection}
\cite{YM-BS:10b}, and \emph{replay attacks} \cite{YM-BS:10a} as
special cases.
Second, we focus on the descriptor model of a power system and we show the
fundamental limitations of static and dynamic detection and identification
procedures. Specifically, we show that static detection procedures are
unable to detect any attack affecting the dynamics, and that attacks
corrupting the measurements can be easily designed to be undetectable. On
the contrary, we show that undetectability in a dynamic setting is much
harder to achieve for an attacker. Specifically, a cyber-physical attack is
undetectable if and only if the attackers' input signal excites uniquely
the zero dynamics of the input/output system. (As a complementary result,
our work \cite{FP-AB-FB:10u} gives necessary and sufficient graph-theoretic
conditions for the absence of zero dynamics, and hence for the absence of
undetectable attacks.) %
Third, we propose a detection and identification procedure based on
geometrically-designed residual filters. Under the assumption of attack
identifiability, our method correctly identifies the attacker set
independently of its strategy. From a system-theoretic perspective, correct
identification is implied by the absence of zero dynamics in our proposed
identification filters. Our design methodology is applicable to linear
systems with direct input to output feedthrough, and it generalizes the
construction presented in \cite{MAM-GCV-ASW-CM:89}.
Fourth and finally, we illustrate the potential impact of our theoretical
results on the standard IEEE 14 bus system (cf. Fig. \ref{fig:ieee14}). For
this system it is known~\cite{YL-MKR-PN:09} that an attack against the
measurement data may remain undetected by a static procedure if the
attacker set compromises as few as four measurements. We show here instead
that such an attack is always detectable by our dynamic detection procedure
provided that at least one bus voltage angle or one generator rotor angle
is measured exactly.

We conclude with two remarks on our contributions.
First, our results (the notions of detectability and identifiability,
the fundamental limitations of static versus dynamic monitoring, and
the geometric design of detection and identification filters) are
analogously and immediately applicable to arbitrary index-one
descriptor systems, thereby including any linear system
$\dot{x}=Ax+Bu$, $y=Cx+Du$, with attack signal $u$. Second, although
we treat here the noiseless case, it is well known~\cite{MB-IVN:93}
that our deterministic detection filters are the key ingredient,
together with Kalman filtering and hypothesis testing, in the design
of statistical identification methods.

\textbf{Organization} Section \ref{sec:model} presents the descriptor
system model of a power network, our framework for the modeling of
cyber-physical attacks, and the detection and identification
problem. Section \ref{sec:static_dynamic} states the fundamental
limitations of static and dynamic detection procedures.  Section
\ref{sec:dynamic_detection} presents the residual filters for dynamic
detection and identification. Section \ref{sec:example} contains the IEEE
14 bus system case study.

%%% Local Variables: 
%%% mode: latex
%%% TeX-master: "main"
%%% End:

\section{Cyber-physical attacks on power networks}
\label{sec:model}
\subsection{Structure-preserving power network model with cyber and
  physical attacks}
We consider the linear small-signal version of the classical
structure-preserving power network model \cite{PWS-MAP:98}. This descriptor
model consists of the {\it dynamic linearized swing equation} and the {\it
  algebraic DC power flow equation}. A detailed derivation from the
nonlinear structure-preserving power network model can be found, for
instance, in \cite{ES:04,FP-AB-FB:10u}.

Consider a connected power network consisting of $n$ generators
$\{g_{1},\dots,g_{n}\}$, their associated $n$ generator terminal buses
$\{b_{1},\dots,b_{n}\}$, and $m$ load buses
$\{b_{n+1},\dots,b_{n+m}\}$. The interconnection structure of the
power network is encoded by a connected admittance-weighted graph. The
generators $g_{i}$ and buses $b_{i}$ form the vertex set of this
graph, and the edges are given by the transmission lines
$\{b_{i},b_{j}\}$ weighted by the susceptance between buses $b_{i}$
and $b_{j}$, as well as the internal connections $\{g_{i},b_{i}\}$
weighted by the transient reactance between each generator $g_{i}$ and
its terminal bus $b_{i}$. The Laplacian matrix associated to the
admittance-weighted graph is the symmetric matrix $\left[
\begin{smallmatrix}
  \subscr{\mc L}{gg} & \subscr{\mc L}{gl}\\
  \subscr{\mc L}{lg} & \subscr{\mc L}{ll}
\end{smallmatrix}
\right] \in \mathbb R^{(2n+m) \times (2n+m)}$, where the first $n$
entries are associated with the generators and the last $n+m$ entries
correspond to the buses. The differential-algebraic model of the power
network is given by the linear continuous-time descriptor\,system
\begin{align}
  \label{eq: power network descriptor system model}
  	E\dot x(t) = A x(t) + P(t), 
\end{align} 
where the state $x=[\delta^\transpose \; \omega^\transpose \;
\theta^\transpose]^\transpose \in \mathbb R^{2n + m}$ consists of the
generator rotor angles $\delta \in \mathbb R^{n}$, the frequencies $\omega
\in \mathbb R^{n}$, and the bus voltage angles $\theta \in \mathbb
R^{m}$. The input term $P(t)$ is due to {\em known} changes in mechanical
input power to the generators or real power demand at the
loads. Furthermore, the descriptor system matrices are
\begin{align}\label{a_matrix}
  E=
  \begin{bmatrix}
         I  & 0  & 0\\
         0 & M & 0\\
         0 & 0  & 0    
  \end{bmatrix}
       ,\,
  A=-
  \begin{bmatrix}
         0 & -I & 0\\
         \subscr{\mc L}{gg} & \subscr{D}{g} & \subscr{\mc L}{gl}\\
         \subscr{\mc L}{lg} & 0 & \subscr{\mc L}{ll}    
  \end{bmatrix},
\end{align}
% matrix is $E=\text{blk-diag}(I,M,0)$,
where $M$ (resp. $D_\textup{g}$) is the diagonal matrix of the
generators' inertias (resp. damping constants).
% with $M$ being the diagonal matrix
% of the generators' inertias, and the system matrix equal to
% \begin{align}\label{a_matrix}
%   A=-\left[\begin{array}{ccc}
%          0 & -I & 0\\
%          \subscr{\mc L}{gg} & \subscr{D}{g} & \subscr{\mc L}{gl}\\
%          \subscr{\mc L}{lg} & 0 & \subscr{\mc L}{ll}\\
%        \end{array}\right]
% %       \in \mathbb R^{(2n+m) \times (2n+m)}
%        \,,
% \end{align} 
% where $\subscr{D}{g} = \textup{diag}(\subscr{D}{g,$1$},\dots,\subscr{D}{g,$n$})$ is the diagonal matrix of
% generator damping constants $\subscr{D}{g,$i$}>0$.
The dynamic and algebraic equations of the linear descriptor system
\eqref{eq: power network descriptor system model} are classically referred
to as the linearized swing equation and the DC power flow equation,
respectively. Notice that the initial condition of system \eqref{eq: power
  network descriptor system model} needs to obey the algebraic constraint
$\subscr{\mc L}{lg} \delta(0) + \subscr{\mc L}{ll} \theta(0) =
P_{\theta}(0)$, where $P_{\theta}(0)$ is the vector containing the entries
$\{2n+1,\dots,2n+m\}$ of $P(0)$.  Finally, we assume the parameters of the
power network descriptor model \eqref{eq: power network descriptor system
  model} to be known, and we remark that they can be either directly
measured, or estimated through dynamic identification techniques, e.g., see
\cite{AC-JHC-AS:11}.

Throughout the paper, the assumption is made that a combination of the
state variables of the descriptor system \eqref{eq: power network
  descriptor system model} is being continuously measured over
time. Let $C \in \mathbb R^{p \times n}$ be the output matrix and let
$y(t) = C x(t)$ denote the $p$-dimensional measurements
vector. Moreover, we allow for the presence of {\em unknown}
disturbances affecting the behavior of the plant \eqref{eq: power
  network descriptor system model}, which, besides reflecting the
genuine failure of network components, can be the effect of a
cyber-physical attack against the network. We classify these
disturbances into \emph{state attacks}, if they show up in the
measurements vector after being integrated through the network
dynamics, and \emph{output attacks}, if they corrupt directly the
measurements vector.\footnote{Because of the linearity of \eqref{eq:
    power network descriptor system model}, the known input $P(t)$ can
  be neglected, since it does not affect the detectability of unknown
  input attacks.}
% We classify these disturbances
% into \emph{physical}, if they are caused by the failure of a network
% component, and \emph{cyber}, if they are caused by the action of a
% \flomargin{think about terminology ``cyber'' and ``physical''}
% malignant attacker against the measurement structure.
% In other words, we make the assumption that physical failures show up
% in the measurements vector after their effects are integrated through
% the dynamics of the system, while cyber attacks directly corrupt the
% system output. 
The network dynamics in the presence of a cyber-physical attack can be
written as
\begin{align}
\label{eq: cyber_physical_fault}
\begin{split}
  E \dot x(t)&=Ax(t) +  
  \underbrace{
    \begin{bmatrix}
%      F_1 & \cdots & F_k \;&\; 0 & \cdots & 0
	F & 0
    \end{bmatrix}
  }_{B}
  \underbrace{
    \begin{bmatrix}
%      f_1 & \ell_{1} \\ \vdots \\ f_k & \ell_{k}
	f(t) \\ \ell(t)
    \end{bmatrix}
    }_{u(t)}
	\,,\\
  y(t)&=Cx(t) + 
  \underbrace{
    \begin{bmatrix}
%      0 & \cdots & 0 \;&\; L_1 & \cdots & L_q
	0 & L
    \end{bmatrix}
  }_{D}
  \underbrace{
    \begin{bmatrix}
%      f_{1} & \ell_1 \\ \vdots \\ f_{k} & \ell_k
	f(t) \\ \ell(t)
    \end{bmatrix}
    }_{u(t)}
	\,.\\
  \end{split}
\end{align}
% \begin{align}\label{eq: cyber_physical_fault}\begin{split}
%   E \dot x(t)&=Ax(t) + B u(t) + \sum_{i=1}^k F_i f_i(t),\\
%   y(t)&=Cx(t) + \sum_{i=1}^q L_i \ell_i (t).
%   \end{split}
% \end{align}
The input signals $f(t)$ and $\ell(t)$ are referred to as {\em state}
and {\em output attack modes}, respectively. The attack modes are
assumed to be unknown and piece-wise continuous functions of time of
dimension $2n+m$ and $p$, respectively, and they act through the full
rank matrices $F \in \mathbb R^{(2n+m) \times (2n+m)}$ and $L \in
\mathbb R^{p \times p}$.
%
% The system \eqref{eq: cyber_physical_fault} is said to be {\em under
%   attack} if there is an active attack mode, i.e., if at least one
% component of $u(t) =
% [f(t)^{\transpose}\,,\,\ell(t)^{\transpose}]^{\transpose}$ is
% non-zero.
For notational convenience, and without affecting generality, we
assume that each state and output variable can be independently
compromised by an attacker. Therefore, we let $F$ and $L$ be the
identity matrices of dimensions $2n+m$ and $p$.
The attack mode $u(t)$ depends upon the specific attack profile. In
the presence of $k \in \mathbb{N}_0$, $k \le 2n+m+p$, attackers
indexed by the {\em attack set} $K \subseteq \until{2n+m+p}$, the
corresponding (vector) attack mode $t \mapsto u_{K}(t) \in \mathbb
R^{2n+m+p}$ has exactly $k$ nonzero entries $u_{K,i}(t)$ for $i \in
K$.
%exactly those $k$ entries $u_{i}(t)$ are nonzero, with $i \in K$, and the
%corresponding (vector) attack mode is termed $u_{K}(t)$. 
Accordingly,\,\,the pair $(B_{K},D_{K})$ is called {\em attack signature}, where $B_K$
and $D_K$ are the submatrices of $B$ and $D$ with columns indexed by
$K$.
% In the absence of attacks, we set $u(t) = 0$ for all $t \in \mathbb
% R_{\geq 0}$.

The model \eqref{eq: cyber_physical_fault} is very general, and it can
capture the occurrence of several concurrent contingencies in the
power network, which are caused either by components failure or
external attacks.\footnote{Genuine failures are a subcase of
  intentional cyber-physical attacks.} For instance,
\begin{enumerate}
\item a change in the mechanical power input to generator $i$
  (resp. in the real power demand of load $j$) is described by the
  attack signature $(B_i,0)$ (resp.  $(B_{2n+j},0)$), and a non-zero
  attack mode $u_{n+i}(t)$ (resp. $u_{2n+j} (t)$);
  % by an appropriate non-zero input signal $f_{n+i}(t)$ with input
  % matrix $B_i$ (resp. $f_{2n+j}(t)$ with $B_j$);
\item a line outage occurring on the line $\{r,s\}$ is modeled by the
  signature $([B_r \; B_s],[0 \; 0])$ and a non-zero mode $[u_r(t) \;
  u_s(t)]^\transpose$ \cite{ES:04}; and
  % $[F_r,F_s]$, and the input signal $f(t) = [k_r(t) \;
  % k_s(t)]^T$, where $k_r(t)$, $k_s(t)$ depend upon the network
  % parameters \cite{ES:04}; and
\item the failure of sensor $i$, or the corruption of the
  $i$-th measurement by an attacker is captured by the signature
  $(0,D_{2n+m+i})$ and a non-zero mode $u_{2n+m+i} (t)$.
\end{enumerate}
% These examples can be readily converted to the corresponding attack
% sets $K$, modes $u_{K}(t)$, and signatures $(B_{K},D_{K})$.

%%% Local Variables: 
%%% mode: latex
%%% TeX-master: "main"
%%% End:

\subsection{Notions of detectability and identifiability for attack sets}
\label{sec:setup}
In this section we present the problem under investigation and we
recall some definitions. Observe that a cyber-physical attack may
remain undetected from the measurements if there exists a normal
operating condition of the network under which the output would be the
same as under the perturbation due to the attacker. Let $x(x_0,u,t)$
denote the network state trajectory generated from the initial state
$x_0$ under the attack signal $u(t)$, and let $y(x_0,u,t)$ be the
output sequence for the same initial condition and input. Throughout
the paper, let $T \subseteq \real_{\ge 0}$ denote the set of time
instants at which the presence of attacks against the network is
checked.

\begin{definition}[\bf Undetectable attack
  set]\label{undetectable_input}
  For the linear descriptor system \eqref{eq: cyber_physical_fault},
  the attack set $K$ is \emph{undetectable} if there exist initial
  conditions $x_1 , x_2 \in \real^{2 n + m}$, and an attack mode $u_K
  (t)$ such that, for all $t \in T$, $y(x_1,u_K,t) =
  y(x_2,0,t)$.
\end{definition}
%Otherwise, the set $K$ is said to be {\em detectable}.
\smallskip

A more general concern than detection is identifiability of attackers,
i.e., the possibility to distinguish from measurements between the
action of two distinct attacks.

\begin{definition}[\bf Unidentifiable attack set]\label{unidentifiable_input}
  For the linear descriptor system \eqref{eq: cyber_physical_fault}, the
  attack set $K$ is \emph{unidentifiable} if there exists an attack set
  $R$, with $|R|\le|K|$ and $R\neq{K}$, initial conditions $x_K , x_R \in
  \real^{2 n + m}$, and attack modes $u_K(t)$, $u_R (t)$ such that, for all
  $t \in T$, $y(x_K,u_{K},t) = y(x_R,u_{R},t)$.
%      For the linear descriptor system \eqref{eq: cyber_physical_fault},
%  the attack set $K$ is \emph{unidentifiable} if there exists a non-identical attack
%  set $R \neq K$, with $|R| \le |K|$, initial conditions $x_K , x_R \in \real^{2 n + m}$, and attack
%  modes $u_K(t)$, $u_R (t)$ such that, for all $t \in \real_{\ge
%    0}$, it holds $y(x_K,u_{K},t) = y(x_R,u_{R},t)$.
  % the attack mode $u_{K}(t)$, with attack set $K$ and attack signature
  % $(B_{K},D_{K})$, is \emph{unidentifiable} if there exist a second
  % attack mode $u_{R}(t)$, with disjoint attack set $R \neq K$ and
  % attack signature $(B_{R},D_{R})$,
  % and initial conditions $x_K , x_R \in \real^{2
  % n + m}$ such that for all $t \in \real_{\ge 0}$ it holds
  % $y(x_K,u_{K},t) = y(x_R,u_{R},t)$. Otherwise, the attack mode $u_K$ is
  % said to be {\em identifiable}.
\end{definition}
\smallskip

Of course, an undetectable attack is also unidentifiable, since it
cannot be distinguished from the zero input. The converse does not
hold. The security problem we consider in this paper is as follows.

\noindent{\em Problem:} % * environment without amsthm package
{(\bf Attack detection and identification)}
For the linear descriptor system \eqref{eq: cyber_physical_fault},
design an attack detection and identification procedure.

Definitions~\ref{undetectable_input} and~\ref{unidentifiable_input} are
immediately applicable to arbitrary constrol systems subjects to external
attacks.  Before proposing a solution to the Attack detection and
identification Problem, we motivate the use of a dynamic detection and
identification algorithm by characterizing the fundamental limitations of
static and dynamic procedures.

\section{Limitations of static and dynamic procedures for detection and
  identification}\label{sec:static_dynamic}
The objective of this section is to show that some fundamental
limitations of a static detection procedure can be overcome by
exploiting the network dynamics. We start by deriving a reduced state
space model for a power network, which is convenient for illustration
and analysis purposes.

\subsection{Kron-reduced representation of a power network}
For the system \eqref{eq: cyber_physical_fault}, let $F = \left[
\begin{smallmatrix}
  F_{\delta}^\transpose &
  F_{\omega}^\transpose &
  F_{\theta}^\transpose
\end{smallmatrix}
\right]^\transpose
$,
$L = 
\left[ 
\begin{smallmatrix}
  L_{\delta}^\transpose & L_{\omega}^\transpose &
  L_{\theta}^\transpose
\end{smallmatrix}
\right]^\transpose $, and
$
C =
\left[
\begin{smallmatrix}
C_\delta & C_\omega & C_\theta
\end{smallmatrix}
\right]$, where the partitioning reflects the state $x =
[\delta^\transpose \; \omega^\transpose \;
\theta^\transpose]^\transpose$. Since the network Laplacian matrix is
irreducible (due to connectivity), the submatrix $\subscr{\mc L}{ll}$
in \eqref{a_matrix} is invertible and the bus voltage angles
$\theta(t)$ can be expressed via the generator rotor angles
$\delta(t)$ and the state attack mode $f(t)$ as
\begin{align}
\label{eq: elimination of bus voltages}
  \theta (t) = - \subscr{\mc L}{ll}^{-1} \subscr{\mc L}{lg} \delta(t) -
  \subscr{\mc L}{ll}^{-1} F_\theta f(t)
 .
\end{align}
Hence, the descriptor system \eqref{eq: cyber_physical_fault} is of
index one \cite{ES:04}. The elimination of the algebraic variables
$\theta(t)$ in the descriptor system \eqref{eq: cyber_physical_fault}
leads to the state space system
\begin{align}
  \label{eq: power network Kron-reduced model}
  \begin{bmatrix}
    \dot \delta \\ \dot \omega
  \end{bmatrix}
  =&\; 
%  \tilde A
\underbrace{
  \begin{bmatrix}
    0 & I \\ - M^{-1} \bigl( \subscr{\mc L}{gg} - \subscr{\mc L}{gl} \subscr{\mc L}{ll}^{-1}
\subscr{\mc L}{lg} \bigr) & -M^{-1}\subscr{D}{g}
  \end{bmatrix}
    }_{\tilde A}
  \begin{bmatrix}
    \delta \\ \omega
  \end{bmatrix}
  \nonumber\\&\;  +
%\tilde B
\underbrace{
\begin{bmatrix}
    F_\delta & 0\\
    M^{-1} F_\omega - M^{-1} \subscr{\mc L}{gl} \subscr{\mc L}{ll}^{-1} F_{\theta} & 0
    \end{bmatrix}
  }_{\tilde B} u
  \,,\\
  y (t) =&\;
%  \tilde C
\underbrace{
  \begin{bmatrix}
    C_\delta  - C_\theta \subscr{\mc L}{ll}^{-1} \subscr{\mc L}{lg} & C_\omega
  \end{bmatrix}
  }_{\tilde C}
  \begin{bmatrix}
    \delta \\ \omega
  \end{bmatrix}
+
%  \tilde D
\underbrace{
  \begin{bmatrix}
    -C_\theta \subscr{\mc L}{ll}^{-1} F_\theta  & L
  \end{bmatrix}
}_{\tilde D} u. \nonumber
\end{align}
This reduction of the passive bus nodes is known as Kron reduction in
the literature on power networks and circuit theory
\cite{FD-FB:11d}. In what follows, we refer to \eqref{eq: power
  network Kron-reduced model} as the {\em Kron-reduced
  system}. Accordingly, for each attack set $K$, the attack signature
$(B_K,D_K)$ is mapped to the corresponding signature $(\tilde B_K,
\tilde D_K)$ in the Kron-reduced system through the transformation for
the matrices $B$ and $D$ described in \eqref{eq: power network
  Kron-reduced model}. Clearly, for any state trajectory of the
Kron-reduced \eqref{eq: power network Kron-reduced model}, the
corresponding state trajectory of the (non-reduced) descriptor power
network model \eqref{eq: cyber_physical_fault} can recovered by
identity \eqref{eq: elimination of bus voltages}.
  
We point out the following subtle but important facts, which are
easily visible in the Kron-reduced system \eqref{eq: elimination of
  bus voltages}. First, a state attack $F_{\theta} f(t)$ on the buses
affects directly the output $y(t)$. Second, for a connected bus
network, the lower block of $\tilde A$ is a fully populated Laplacian
matrix, and $\subscr{\mc L}{ll}^{-1}$ and $\subscr{\mc L}{gl}
\subscr{\mc L}{ll}^{-1}$ are both positive matrices
\cite{FD-FB:11d}. As one consequence, an attack on a single bus
affects the {\em entire} network and not only the locally attacked
node or its vicinity. Third and finally, the mapping from the input
signal $u(t)$ and the initial condition $x(0)$ (subject to the
constraint \eqref{eq: elimination of bus voltages} evaluated at $t =
0$) to the output signal $y(t)$ of the descriptor system \eqref{eq:
  cyber_physical_fault} coincides with the corresponding input and
initial state to output map of the associated Kron-reduced system
\eqref{eq: power network Kron-reduced model}. Hence, the definition of
identifiability (resp. detectability) of an attack set is analogous
for the Kron-reduced system \eqref{eq: power network Kron-reduced
  model}, and we can directly state the following lemma.

%\begin{lemma}{\bf(Equivalence of identifiable attack):}
\begin{lemma}{\bf(Equivalence of detectability and identifiability under Kron reduction):}
\label{lemma:equivalence}
For the power network descriptor system \eqref{eq:
  cyber_physical_fault}, the attack set $K$ is identifiable
(resp. detectable) if and only if it is identifiable
(resp. detectable) for the associated Kron-reduced system \eqref{eq:
  power network Kron-reduced model}.
\end{lemma}
%\begin{proof}
%%  The lemma follows from the fact that the input and initial
%%  condition to output map of the system \eqref{eq:
%%    cyber_physical_fault} coincides with the input to output map of
%%  the associated Kron-reduced system and equation \eqref{eq:
%%    elimination of bus voltages} at time $t = 0$.
%    The lemma follows from the fact that the mapping from the input
%    signal $u(t)$ and the initial condition $x(0)$ (subject to the
%    constraint \eqref{eq: elimination of bus voltages} evaluated at
%    $t = 0$) to the output signal $y(t)$ of the descriptor system
%    \eqref{eq: cyber_physical_fault} coincides with the corresponding
%    input and initial state to output map of the associated
%    Kron-reduced system \eqref{eq: power network Kron-reduced model}.
%\end{proof}

Following Lemma \ref{lemma:equivalence}, we study detectability and
identifiability of attacks against the power network descriptor model
\eqref{eq: cyber_physical_fault} by analyzing the associated
Kron-reduced system \eqref{eq: power network Kron-reduced model}.

\subsection{Fundamental limitations of a Static Detector}
By \emph{Static Detector}, or, with the terminology of
\cite{AA-AGS:04}, \emph{Bad Data Detector}, we denote an algorithm
that uses the network measurements to check for the presence of
attacks at some predefined instants of time, and without exploiting
any relation between measurements taken at different time instants. By
Definition \ref{undetectable_input}, an attack is undetectable by a
Static Detector if and only if, for all time instances $t$ in a
countable set $T$, there exists a vector $\xi(t)$ such that $y(t) =
\tilde C \xi(t)$. Without loss of generality, we set $T = \mathbb
N$. Loosely speaking, the Static Detector checks whether, at a
particular time instance $t \in \mathbb N$, the measured data is
consistent with the measurement equation, for example, the power flow
equation at a bus. Notice that our definition of Static Detector is
compatible with \cite{YL-MKR-PN:09}, where an attack is detected if
and only if the residual $r(t) = y(t) - \tilde C
[\hat\delta(t)^{\transpose} \;
\hat\omega(t)^{\transpose}]^{\transpose}$ is nonzero for some $t \in
\mathbb N$, where $[\hat\delta(t)^{\transpose} \;
\hat\omega(t)^{\transpose}]^{\transpose} = \tilde{C}^{\dag}
y(t)$. % is
% the best unbiased estimate of the network state at time $t$.
% As pointed out in \cite{YL-MKR-PN:09},
If $r(t) \neq 0$, then a malfunction is detected, and it is undetected
otherwise.\footnote{Similar conclusion can be drawn for the case of
  noisy measurements.}

% In particular, after computing an estimate
% $[\hat\delta(t)^{\transpose} \;
% \hat\omega(t)^{\transpose}]^{\transpose} = \tilde{C}^{\dag} y(t)$ of
% the network state at time $t$, the Static Detector analyzes the
% residual $r(t) = y(t) - \tilde C [\hat\delta(t)^{\transpose} \;
% \hat\omega(t)^{\transpose}]^{\transpose}$.

\begin{theorem}
\emph{(\bf Static detectability of cyber-physical attacks)}
\label{Theorem: Static detectability of cyber-physical attacks}
For the power network descriptor system \eqref{eq:
  cyber_physical_fault} and an attack set $K$, the following two statements are
equivalent:
\begin{enumerate}
  \item the attack set $K$ is undetectable by a Static Detector;
  \item there exists an attack mode $u_K (t)$ such that, for some
    $\delta (t)$ and $\omega(t)$, at every $t \in \mathbb N$ it holds
    \begin{align}\label{eq:cond2_static}
      \tilde C 
      \begin{bmatrix}
        \delta(t)\\
        \omega(t)
      \end{bmatrix} + \tilde D u_K(t) = 0\,,
    \end{align}
    where $\tilde C$ and $\tilde D$ are as in \eqref{eq: power network
      Kron-reduced model}.
  % \item there exists an output attack mode $\ell_K(t)$\fpmargin{define
  %     $f_K$ and $\ell_K$} such that $L \ell_K(t_i) \in \Image (C)$ at
  %   every instant $t_i \in T$.
\end{enumerate}
Moreover, there exists an attack set $K$ undetectable by a Static Detector if and only if there exist $x \in \real^{2n}$ and $g \in \real^{|K|}$ such that $\tilde C x + \tilde D_K g =
0$.
%Moreover, there exists a static-undetectable attack set $K$, with $|K|
%= k \in \mathbb{N}$, if and only if there exist $g \in \real^{|K|}$
%and $0 \neq x \in \real^{2n}$ such that $\tilde C x + \tilde D_K g =
%0$.
%
% a necessary and sufficient condition for the existence of an
% undetectable attack set $K$ is that, for some $g \in \real^{|K|}$ and
% $0 \neq x \in \real^{2n}$, it holds $\tilde C x + \tilde D_K g = 0$.
% let $K$ be the attack set, and let $u_K (t)$
% be the attack mode. The following three statements are equivalent:
% \begin{enumerate}
% \item the cyber-physical attack mode $u_K(t)$ is undetectable by a
%   Static Detector;
% \item\label{static2} for each $t_i \in T$, there exist $\delta (t_i)$ and
%   $\omega(t_i)$ such that
%   \begin{align*}
%     \tilde C 
%     \begin{bmatrix}
%       \delta(t_i)\\
%       \omega(t_i)
%     \end{bmatrix} + \tilde D u_K(t_i) = 0\,,
%   \end{align*}
%   where $\tilde C$ and $\tilde D$ are as in \eqref{eq: power network
%     Kron-reduced model};
% \item\label{static3} for each $t_i \in T$, it holds $L \ell(t_i)
%   \in \Image (C)$.
% \end{enumerate}
% Moreover, a necessary and sufficient condition for the existence of an
% undetectable attack is that, for some $g \in \real^{2n+m+p}$ and some nonzero $x \in
% \real^{2n}$, it holds
% $\tilde C x + \tilde D g = 0$.
\end{theorem}
\smallskip

Before presenting a proof of the above theorem, we highlight that a
necessary and sufficient condition for the equation
\eqref{eq:cond2_static} to be satisfied is that $L \ell(t) \in \Image
(C)$ at all times $t \in \mathbb N$, where $\ell(t)$ is the output
attack mode, i.e., the vector of the last $p$ components of
$u_K(t)$. Hence, statement (ii) in Theorem \ref{Theorem: Static
  detectability of cyber-physical attacks} implies that {\em no} state
attack can be detected by a static detection procedure, and that an
undetectable output attack exists if and only if $\Image (D_{K}) \cap
\Image (C) \neq \{ 0 \}$.

\begin{pfof}{Theorem \ref{Theorem: Static detectability of
      cyber-physical attacks}}
  % Due to invariance of detectability under Kron reduction by Lemma
  % \ref{lemma:equivalence}, Theorem \ref{Theorem: Static detectability
  %   of cyber-physical attacks} will be proved for the Kron-reduced
  % system \eqref{eq: power network Kron-reduced model} instead of the
  % descriptor system \eqref{eq: cyber_physical_fault}.
 %
  As previously discussed, the attack $K$ is undetectable by a Static
  Detector if and only if for each $t \in \mathbb N$ there exists
  $\delta(t)$, $\omega(t)$, and $u_{K}(t)$ such that
  \begin{align*}
    r(t) = y(t) - \tilde C
    \tilde C^{\dagger}
    y(t) = (I-\tilde C \tilde C^{\dagger}) \left( \tilde C
    \begin{bmatrix}
      \delta(t) \\ \omega(t)
    \end{bmatrix}
    + \tilde D u_K(t) \right)
  \end{align*}
  vanishes. Consequently, $r(t) = (I - \tilde C \tilde C^{\dagger})
  \tilde D u_K(t)$, and the attack set $K$ is undetectable
  if and only if $\tilde Du_{K}(t) \in \Image(\tilde C)$, which is
  equivalent to statement (ii). The last necessary and sufficient condition
  in the theorem follows from (ii).
\end{pfof}

We now focus on the static identification problem. Following Definition
\ref{unidentifiable_input}, the following result can be asserted.

\begin{theorem}\emph{(\bf Static identification of cyber-physical attacks)}
  \label{Theorem: Static identifiability of cyber-physical attacks}
  For the power network descriptor system \eqref{eq: cyber_physical_fault}
  and an attack set $K$, the following two statements are equivalent:
  \begin{enumerate}
  \item the attack set $K$ is unidentifiable by a Static Detector;
  \item there exists an attack set $R$, with $|R|\le|K|$ and $R\neq{K}$,
    and attack modes $u_K (t)$, $u_R (t)$, such that, for some $\delta (t)$
    and $\omega(t)$, at every $t \in \mathbb N$, it holds
    \begin{align*}
      \tilde C 
      \begin{bmatrix}
        \delta(t)\\
        \omega(t)
      \end{bmatrix} + \tilde D \left( u_K(t) + u_R(t) \right)= 0\,,
    \end{align*}
    where $\tilde C$ and $\tilde D$ are as in \eqref{eq: power network
      Kron-reduced model}.
  \end{enumerate}
  Moreover, there exists an attack set $K$ unidentifiable by a Static
  Detector if and only if there exists an attack set $\bar K$, $|\bar
  K| \le 2 |K|$, which is undetectable by a Static Detector.
  % if and only if there exists a set $\bar K$, with $|\bar K| \le 2k$,
  % $g \in \real^{|\bar K|}$ and $0 \neq x \in \real^{2n}$, it holds $\tilde C
  % x + \tilde D_{\bar K} g = 0$.
\end{theorem}
\smallskip

Similar to the fundamental limitations of static detectability in
Theorem \ref{Theorem: Static detectability of cyber-physical attacks},
Theorem \ref{Theorem: Static identifiability of cyber-physical
  attacks} implies that, for instance, state attacks cannot be
identified and that an undetectable output attack exists if and only
if $\Image (D_{\bar K}) \cap \Image (C) \neq \{ 0 \}$.
 
\begin{pfof}{Theorem \ref{Theorem: Static identifiability of cyber-physical attacks}}
  Because of the linearity of the system \eqref{eq:
    cyber_physical_fault}, the unidentifiability condition in
  Definition \ref{unidentifiable_input} is equivalent to
  $y(x_K-x_R,u_{K}-u_{R},t) = 0$, for some initial condition $x_K$,
  $x_R$, and attack mode $u_K (t)$, $u_R(t)$. The equivalence between
  statements (i) and (ii) follow. The last statement follows from
  Theorem \ref{Theorem: Static detectability of cyber-physical
    attacks}.
\end{pfof}

\subsection{Fundamental limitations of a Dynamic Detector}
In the following we refer to a security system having access to the
{\em continuous} time measurements signal $y(t)$, $t \in \real_{\ge
  0}$, as a {\em Dynamic Detector}. As opposed to a Static Detector, a
Dynamic Detector checks for the presence of attacks at every instant
of time $t \in \real_{\geq 0}$. By Definition
\ref{undetectable_input}, an attack is undetectable by a Dynamic
Detector if and only if there exists a network initial state $\xi(0)
\in \real^{2n}$ such that $y(t) = \tilde C e^{\tilde A t}\xi(0)$ for
all time instances $t\in \real_{\ge 0}$. Intuitively, a Dynamic
Detector is harder to mislead than a Static Detector.
%The following theorem formalizes this expected result.

\begin{theorem}
\emph{(\bf Dynamic detectability of cyber-physical attacks)}
\label{Theorem: Dynamic detectability of cyber-physical attacks}
For the power network descriptor system \eqref{eq:
  cyber_physical_fault} and an attack set $K$, the following two
statements are equivalent:
\begin{enumerate}
  \item the attack set $K$ is undetectable by a Dynamic Detector;
  \item there exists an attack mode $u_K (t)$ such that, for some
    $\delta (0)$ and $\omega(0)$, at every $t \in \real_{\ge 0}$, it
    holds
        \begin{align*}
      \tilde C 
      e^{\tilde A t}
      \begin{bmatrix}
        \delta(0)\\
        \omega(0)
    \end{bmatrix}
    \! +  \tilde C \! \int_0^t \! e^{\tilde A (t-\tau)} \tilde B u_K (\tau)
    d\tau 
  =- \tilde D u_K(t),
  \end{align*}
where $\tilde A$, $\tilde B$, $\tilde C$, and $\tilde D$ are as in
  \eqref{eq: power network Kron-reduced model}.
\end{enumerate}
Moreover, there exists an attack set $K$ undetectable by a Dynamic
Detector if and only if there exist $s \in \complex$, $g \in
\real^{|K|}$, and $x \in \real^{2n}$, $x \neq 0$, such that
$(sI-\tilde A)x - \tilde B_K g = 0$ and $\tilde C x + \tilde D_K g =
0$.
%  a necessary and sufficient condition for the existence of an
% undetectable attack set $K$ is that, for some $s \in \complex$, $g \in
% \real^{|K|}$, and some nonzero $x \in \real^{2n}$, it holds
% $(sI-\tilde A)x - \tilde B g = 0$ and $\tilde C x + \tilde D_K g = 0$.
% let $K$ be the attack set, and let $u_K(t)$
% be the attack mode. The following two statements are equivalent:
% \begin{enumerate}
% \item the cyber-physical attack mode $u_K(t)$ is undetectable by a
%   Dynamic Detector;
% \item\label{dynamic2} there exist $\delta (0)$ and $\omega(0)$ such
%   that, for each $t \in \real_{\ge 0}$, it holds
%     \begin{align*}
%       \tilde C 
%       e^{\tilde A t}
%       \begin{bmatrix}
%         \delta(0)\\
%         \omega(0)
%     \end{bmatrix}
%     \! +  \tilde C \! \int_0^t \! e^{\tilde A (t-\tau)} \tilde B u_K (\tau)
%     d\tau 
%   =- \tilde D u_K(t),
%   \end{align*}
%   where $\tilde A$, $\tilde B$, $\tilde C$, and $\tilde D$ are as in
%   \eqref{eq: power network Kron-reduced model}.
% \end{enumerate}
% A \emph{necessary} condition for the input $u_K(t)$ to be
% undetectable\,is
% \begin{enumerate}
%   \item[(iii)] $L \ell(t) \in \Image (C)$ for all $t \in
%   \real_{\ge 0}$.
% \end{enumerate}
% Moreover, a necessary and sufficient condition for the existence of an
% undetectable attack is that, for some $s \in \complex$, $g \in
% \real^{2n+m+p}$, and some nonzero $x \in \real^{2n}$, it holds
% $(sI-\tilde A)x - \tilde B g = 0$ and $\tilde C x + \tilde D g = 0$.
\end{theorem}
\smallskip

Before proving Theorem \ref{Theorem: Dynamic detectability of
  cyber-physical attacks}, some comments are in order. First, state
attacks {\em can be detected} in the dynamic case. Second, an attacker
needs to inject a signal which is consistent with the network dynamics
at every instant of time to mislead a Dynamic Detector. Hence, as
opposed to the static case, the condition $L \ell(t) \in \Image (C)$
needs to be satisfied for every $t \in \real_{\ge 0}$, and it is only
necessary for the undetectability of an output attack. Indeed, for
instance, state attacks can be detected even though they automatically
satisfy the condition $0 = L \ell(t) \in \Image (C)$. Third and
finally, according to the last statement of Theorem \ref{Theorem:
  Dynamic detectability of cyber-physical attacks}, the existence of
invariant zeros for the Kron-reduced system $(\tilde A,\tilde B_K,
\tilde C, \tilde D_K)$ is equivalent to the existence of an
undetectable attack mode $u_K (t)$.\footnote{For the system $(\tilde
  A,\tilde B_K, \tilde C, \tilde D_K)$, the value $s \in \complex$ is
  an invariant zero if there exists $x \in \real^{2n}$, with $x \neq
  0$, $g \in \real^{|K|}$, such that $(sI-\tilde A)x - \tilde B_K g =
  0$ and $\tilde C x + \tilde D_K g = 0$. For a linear dynamical
  system, the existence of invariant zeros is equivalent to the
  existence of zero dynamics \cite{HLT-AS-MH:01}.}
As a consequence, for the absence of undetectable cyber-physical
attacks, a dynamic detector performs better that a static detector,
while requiring, possibly, fewer measurements. A related example is in
Section \ref{sec:example}.

% Conversely, any attack $u_K (t)$ is detectable if and only if the
% Kron-reduced system $(\tilde A,\tilde B, \tilde C, \tilde D)$ has no
% invariant zeros. fewer measurements because C not full column rank,
% necessarily
% This condition suggests that, generally, fewer measurements are needed
% in a dynamic detection procedure than in a static detection
% procedure. In Section \ref{sec:dynamic_detection} we suggest dynamic
% detection and identification procedures that confirm this intuition
% for the IEEE 14 bus system, see Section \ref{sec:example}.

\begin{pfof}{Theorem \ref{Theorem: Dynamic detectability of
    cyber-physical attacks}}
% The theorem will be proved for the Kron-reduced system \eqref{eq:
%   power network Kron-reduced model} instead of the descriptor system
% \eqref{eq: cyber_physical_fault}, since attack detectability is
% invariant under Kron reduction by Lemma \ref{lemma:equivalence}.
By Definition \ref{undetectable_input} and linearity of the system
\eqref{eq: power network Kron-reduced model}, the attack mode $u_K(t)$
is undetectable by a Dynamic Detector if and only if there exists
$[\delta(0)^{\transpose} \; \omega(0)^{\transpose}]^{\transpose}$ such
that $y([\delta(0)^{\transpose} \;
\omega(0)^{\transpose}]^{\transpose},u_K,t) = 0$ for all $t \in
\mathbb R_{\geq 0}$. Hence, statements (i) and (ii) are equivalent.
Following condition (ii) in Theorem \ref{Theorem: Dynamic
  detectability of cyber-physical attacks}, an attack $u_K(t)$ may
remain undetected to a Dynamic Detector if and only if 
$u_K(t)$ is an input-zero for some initial condition.
%\cite{HLT-AS-MH:01}. 
%Recall that, for the system \eqref{eq: power
%  network Kron-reduced model}, there exists an input-zero if and only
%if there exist $s \in \mathbb C$, $x \neq 0$, and $g$ such that
%$(sI-\tilde A)x - \tilde B g = 0$ and $\tilde C x + \tilde D g = 0$.
%%
%This concludes the proof.
\end{pfof}

% \subsection{Fundamental limitations of dynamic identification
%   procedures}
% In the presence of $k \in \until{2n+m+p}$ attackers with attack set $K
% = \{i_1,\dots,i_k\}$ and failure signatures $(B_{K},D_{K})$, we let
% $(\tilde B_{K},\tilde D_{K})$ be the corresponding failure modes in
% the Kron-reduced system \eqref{eq: power network Kron-reduced model}.

We now focus on the identification problem.

\begin{theorem}
\emph{(\bf Dynamic identifiability of cyber-physical attacks)}
\label{Theorem: Dynamic identifiability of cyber-physical attacks}
For the power network descriptor system \eqref{eq:
  cyber_physical_fault}, the following two statements are equivalent:
\begin{enumerate}
\item the attack set $K$ is unidentifiable by a Dynamic Detector;
\item there exists an attack set $R$, with $|R| \le |K|$ and $R\neq{K}$ if
  $|R| = |K|$, and attack modes $u_K(t)$, $u_R(t)$, such that, for some
  $\delta (0)$ and $\omega(0)$, at every $t \in \real_{\ge 0}$, it holds
  \begin{align*}
    &\tilde C 
    e^{\tilde A t}
    \begin{bmatrix}
      \delta(0)\\
      \omega(0)
    \end{bmatrix}
    \! + \tilde C \! \int_0^t \! e^{\tilde A (t-\tau)} \tilde B \left( u_K
    (\tau) + u_R(\tau) \right) d\tau
    \\
    &= - \tilde D \left(u_K(t) + u_R(t) \right)
  \end{align*}
  where $\tilde A$, $\tilde B$, $\tilde C$, and $\tilde D$ are as in
  \eqref{eq: power network Kron-reduced model}.
% \item there exists an attack set $R$, with $|K| = |R|$ and $K \neq
%     R$, such that
%     \begin{align*}
%       \begin{bmatrix}
% 		sI - \tilde A & - \tilde B_{K} & - \tilde B_{R} \\
% 		\tilde C & \tilde D_{k} & \tilde D_{R}
%               \end{bmatrix}
%               \begin{bmatrix}
% 		x \\ g_{K} \\ g_{R}
%               \end{bmatrix}
%               =
%               \begin{bmatrix}
% 		0 \\ 0 \\ 0
% 	\end{bmatrix}
% 	\,,
%       \end{align*}
%       for some $s \in \complex$, $0 \neq x \in \real^{2n}$, $g_K,
%       g_R \in \real^{2n+m+p}$.
\end{enumerate}
Moreover, there exists an attack set $K$ unidentifiable by a Dynamic
Detector if and only if there exists an attack set $\bar K$, 
$|\bar K| \le 2 |K|$, which is unidentifiable by a Dynamic Detector.
% under the cyber-physical attack $u(t)$, $t \in
% \real_{\ge 0}$. 
%   \begin{enumerate}
%
%   \item there exists a cyber-physical attack $u(t)$ of maximum
%     cardinality $k \in \until{2n+m+p}$ which is unidentifiable by a
%     Dynamic Detector;
%   
%   \item there exists two disjoint attack sets $K$ and $R$, each of
%     maximum cardinality $k$, such that for some $s \in \complex$, $x
%     \neq 0$, $g_{K} \in \mathbb R^{K \times 1}$, and $g_{R} \in
%     \mathbb R^{R \times 1}$ it holds that
% \begin{equation*}
% 	\begin{bmatrix}
% 		sI - \tilde A & - \tilde B_{K} & - \tilde B_{R} \\
% 		\tilde C & \tilde D_{k} & \tilde D_{R}
% 	\end{bmatrix}
% 	\begin{bmatrix}
% 		x \\ g_{K} \\ g_{R}
% 	\end{bmatrix}
% 	=
% 	\begin{bmatrix}
% 		0 \\ 0 \\ 0
% 	\end{bmatrix}
% 	\,.
% \end{equation*}
% \end{enumerate}
\end{theorem}

\begin{pf}
  Notice that, because of the linearity of the system \eqref{eq:
    cyber_physical_fault}, the unidentifiability condition in
  Definition \ref{unidentifiable_input} is equivalent to the condition
  $y(x_K-x_R,u_{K}-u_{R},t) = 0$, for some initial condition $x_K$,
  $x_R$, and attack mode $u_K (t)$, $u_R$. The equivalence
  between statements (i) and (ii) follows.
  % the attack mode $u_{K}(t)$ with attack set $K$ (of maximum
  % cardinality $k$) is unidentifiable if and only if there exists a
  % second disjoint attack set $R$ (of maximum cardinality $k$) with
  % the property
  % $\Image([B_{K}^\transpose,D_{K}^\transpose]^\transpose) \cap
  % \Image([B_{R}^\transpose,D_{R}^\transpose]^\transpose) = \{ 0 \}$
  % (equivalent to $R$ and $K$ being disjoint),
  % with attack mode $u_{R}(t)$, and initial conditions $x_1 , x_2 \in
  % \real^{2 n + m}$ such that for all $t \in \real_{\ge 0}$ it holds
  % $y(x_1,u_{K},t) = y(x_2,u_{R},t)$. Due to linearity of the system,
  % the latter condition is equivalent to a zero output
  % $y(x_1-x_2,u_{K}-u_{R},t) = 0$ of the system with input matrix
  % $[B_{K},B_{R}]$ and direct feedthrough matrix $[D_{K},D_{R}]$. The
  % latter condition is again equivalent to statement (ii).
\end{pf}

In other words, the existence of an unidentifiable attack set of
cardinality $k$ is equivalent to the existence of invariant zeros for
the system $(\tilde A,\tilde B_{\bar K}, \tilde C, \tilde D_{\bar
  K})$, for some attack set $\bar K$ with $|\bar K| \le 2k$. A careful
reader may notice that condition (ii) in Theorem 3.4 is hard to verify
because of its combinatorial complexity: % in order to check if the
% attack set $K$ is identifiable,
one needs to certify the absence of invariant zeros for all possible
distinct pairs of $|K|$-dimensional attack sets. Then, a conservative
verification of condition (ii) requires $\binom{2n+m+p}{2|K|}$
tests. In \cite{FP-AB-FB:10u} we partially address this complexity
problem by presenting an intuitive and easy to check graph-theoretic
condition for a given network topology and generic system parameters.

% \begin{remark}(\emph{\bf Stealth, false-data injection, and replay
%     attacks}) The following prototypical attacks can be modeled and
%   analyzed through our theoretical framework. First, stealth attacks
%   as defined in \cite{DG-HS:10} corresponds to output attacks $D_K
%   u_K(t) \in \Image(C)$. Second, (dynamic) false-data injection
%   attacks are presented in \cite{YM-BS:10b}. These attacks are
%   undetectable output attacks rendering the unstable modes (if any) of
%   the system unobservable. Clearly, these unobservable modes are
%   included in the invariant zeros set. Third and finally, for a replay
%   attack \cite{YM-BS:10a}, the attack signature is such that
%   $\Image(C) \subseteq \Image(D_K)$ and $B_K \neq 0$. Notice that the
%   resulting system would have an infinite number of invariant
%   zeros\fpmargin{check this. the number can be finite}, and, hence,
%   replay attacks are only particular attack modes. \oprocend
% \end{remark}

\begin{remark}(\emph{\bf Stealth, false-data injection, and replay
    attacks}) The following prototypical attacks can be modeled and
  analyzed through our theoretical framework:
  \begin{enumerate}
  \item stealth attacks, as defined in \cite{DG-HS:10}, correspond to
    output attacks satisfying $D_K u_K(t) \in \Image(C)$;
  \item (dynamic) false-data injection attacks, as defined in
    \cite{YM-BS:10b}, are output attacks rendering the unstable modes
    (if any) of the system unobservable. These unobservable modes are
    included in the invariant zeros set; and
  \item replay attacks, as defined in \cite{YM-BS:10a}, are state and
    output attacks satisfying $\Image(C) \subseteq \Image(D_K)$, $B_K
    \neq 0$. The resulting system may have an infinite number of
    invariant zeros: if the attacker knows the system model, then it
    can cast very powerful undetectable attacks.
  \end{enumerate}
  In \cite{YM-BS:10a}, a monitoring signal (unknown to the attacker)
  is injected into the system to detect replay attacks. It can be
  shown that, if the attacker knows the system model, and if the
  attack signal enters additively as in \eqref{eq:
    cyber_physical_fault}, then the attacker can design undetectable
  attacks without knowing the monitoring signal. Therefore, 
  % for this
  % attack model
  the fundamental limitations presented in Section
  \ref{sec:static_dynamic} are also valid for \emph{active} detectors,
  which are allowed to inject monitoring signals to reveal
  attacks.\oprocend
\end{remark}

% We remark that condition (ii) in Theorem \ref{Theorem: Dynamic
%   identifiability of cyber-physical attacks} is not constructive and
% hard to verify due to combinatorial complexity. The results in
% \cite{FP-AB-FB:10u} present intuitive graph-theoretic conditions when
% condition (ii) is satisfied for a given network topology and generic
% system parameters.

%%% Local Variables: 
%%% mode: latex
%%% TeX-master: "main"
%%% End:

\section{Design of dynamic detection and identification
  procedures}\label{sec:dynamic_detection} 
%\input{dynamic_detection}
% We now design the filters which constitute the basis of our dynamic
% attack detection and identification procedure.

\subsection{Detection of attacks} 
% Recall from Theorem \ref{Theorem: Dynamic detectability of
%   cyber-physical attacks}, that an attack $u$ is {\em detectable} if
% and only if the Kron-reduced system $(\tilde A,\tilde B, \tilde C,
% \tilde D)$ has no invariant zeros. In this case, we propose the
% following residual filter to successfully identify the attack.
We start by considering the attack detection problem, whose
solvability condition is in Theorem \ref{Theorem: Dynamic
  detectability of cyber-physical attacks}.
% Consider the power network descriptor system \eqref{eq:
% cyber_physical_fault}, and assume that the attack set $K$ is
% detectable (see Definition \ref{undetectable_input}).
We propose the following residual filter to detect cyber-physical
attacks.

\begin{theorem}[\bf Attack detection filter]\label{Proposition:
    Attack detection filter}
  Consider the power network descriptor system \eqref{eq:
    cyber_physical_fault} and the associated Kron-reduced system
  \eqref{eq: power network Kron-reduced model}. Assume that the attack
  set is detectable and that the network initial state $x(0)$ is
  known.  Consider the {\em detection filter} \begin{align}
	\begin{split}
		\dot w(t) &= (\tilde A + G \tilde C) w(t) - Gy(t) , \\
		r(t) &= \tilde C w(t) - y(t) ,
	\end{split}
	\label{eq: detection filter}
\end{align}
where $w(0) = x(0)$, and $G \in \mathbb R^{2n \times p}$ is such that
$\tilde A + G \tilde C$ is a Hurwitz matrix. Then $r (t) = 0$ at all
times $t \in \real_{\ge 0}$ if and only if $u(t) = 0$ at all times $t
\in \real_{\ge 0}$.
% Assume that all
% cyber-physical attacks $u(t)$, $t \in \real_{\ge 0}$ are detectable,
% or equivalently, Kron-reduced model \eqref{eq: power network
%   Kron-reduced model} has no invariant zeros. Consider the {\em
%   detection filter}
% \begin{align}
% 	\begin{split}
% 		\dot w(t) &= (\tilde A + G \tilde C) w(t) - Gy(t) , \\
% 		r(t) &= \tilde C w(t) - y(t) ,
% 	\end{split}
% 	\label{eq: detection filter}
% \end{align}
% where $G \in \mathbb R^{2n \times q}$ is such that $\tilde A + G
% \tilde C$ is a Hurwitz matrix, and set $w(0) = x(0)$. Then the filter
% output $r(t)$ converges exponentially to zero if and only if $u(t) =
% 0$ for all $t \in \real_{\ge 0}$.
\end{theorem}
\begin{pf}
  Consider the error $e(t) = w(t)-x(t)$ between the states of the
  filter \eqref{eq: detection filter} and the Kron-reduced system
  \eqref{eq: power network Kron-reduced model}. The error dynamics
  with output $r(t)$ are then
  \begin{align}
    \begin{split}
      \dot e(t) &= (\tilde A + G \tilde C) e(t) - (\tilde B + G\tilde D )u(t) , \\
      r(t) &= \tilde C e(t) - \tilde D u(t) ,
    \end{split}
    \label{eq: detection filter - error dynamics}
  \end{align}
  where $e(0) = 0$. Clearly, if the error system \eqref{eq: detection
    filter - error dynamics} has no invariant zeros, then $r(t)=0$ for all
  $t \in \real_{\ge 0}$ if and only if $u(t) = 0$ for all $t \in \real_{\ge
    0}$ and the claimed statement is true.  The error system \eqref{eq:
    detection filter - error dynamics} has no invariant zeros if and only
  if there exists no triple $(s,\bar w,g) \in \mathbb C \times \mathbb
  R^{2n} \times \mathbb R^{p}$ satisfying
  \begin{equation}
    \begin{bmatrix}
      sI - (\tilde A+G \tilde C) & \tilde B + G \tilde D  \\
      \tilde C & -\tilde D
    \end{bmatrix}
    \begin{bmatrix}
      \bar w \\ g
    \end{bmatrix}
    =
    \begin{bmatrix}
      0 \\ 0
    \end{bmatrix}
    \label{eq: pencil for error system}
	\,.
      \end{equation}
      The second equation of \eqref{eq: pencil for error system} yields
      $\tilde C x = \tilde D g$. Thus, by substituting $\tilde C x$ by
      $\tilde D g$ in the first equation of \eqref{eq: detection filter -
        error dynamics}, the set of equations \eqref{eq: pencil for error
        system} can be equivalently written as
\begin{equation*}
  \begin{bmatrix}
    sI - \tilde A & \tilde B  \\
	\tilde C & -\tilde D
      \end{bmatrix}
      \begin{bmatrix}
	\bar w \\ g
	\end{bmatrix}
	=
	\begin{bmatrix}
          0 \\ 0
	\end{bmatrix}
	\,.
      \end{equation*}
      Finally, note that the solution $(s,-\bar w,g)$ to the above set of
      equations yields an invariant zero, zero state, and zero input for
      the Kron-reduced system \eqref{eq: power network Kron-reduced
        model}. By the detectability assumption, the Kron-reduced system
      \eqref{eq: power network Kron-reduced model} has no zero dynamics. We
      conclude that the error system \eqref{eq: detection filter - error
        dynamics} has no zero dynamics, and the statement is true.
\end{pf}

In summary, the implementation of the residual filter~\eqref{eq: detection
  filter} guarantees the detection of any detectable attack set.

\subsection{Identification of attacks} 

We now focus on the attack identification problem, whose solvability
condition is in Theorem \ref{Theorem: Dynamic identifiability of
  cyber-physical attacks}. Unlike the detection case, the
identification of the attack set $K$ requires a combinatorial
procedure, since, a priori, $K$ is one of the $\binom{2n+m+p}{|K|}$
possible attack sets. As key component of our identification
procedure, we propose a residual filter to determine whether a
predefined set coincides with the attack set.
% Recall from Theorem \ref{Theorem: Dynamic detectability of
% cyber-physical attacks} that the attack mode $u_{K}(t)$ is
% identifiable if for any disjoint attack set $R \neq K$ the joint
% system $(A,[B_{K} \; B_{R}],C,[D_{K} \; D_{R}])$ has only trivial
% zero dynamics. In what follows, we present an attack identification
% filter.
% For convenience, the filter matrices are defined in an appropriate
% coordinate frame, and, for the implementation in the original
% coordinates, the following definitions are propaedeutic. 

We next introduce in a coordinate-free geometric way the key elements
of this residual filter based on the notion of condition-invariant
subspaces \cite{HLT-AS-MH:01}. Let $K$ be a $k$-dimensional attack
set, and let $\tilde B_K$, $\tilde D_K$ be as defined right after the
Kron reduced model~\eqref{eq: power network Kron-reduced model}. Let
$[V_{K}^\transpose \; Q_{K}^\transpose]^\transpose \in \mathbb R^{p
  \times p}$ be an orthonormal matrix such that
\begin{align*}
  V_{K} = \Basis(\Image(\tilde D_K)), \text{ and } Q_{K} =
  \Basis(\Image(\tilde D_K)^\perp),
\end{align*}
and let 
\begin{align}
  \label{eq: Bz Bk}
  B_{Z} = \tilde B_{K} (V_{K} \tilde D_K)^\dag, \text{ and }\bar B_{K}
  = \tilde B_{K}(I-D_{K}D_{K}^{\dag}).
\end{align}
Define the subspace $\Star \subseteq \real^{2n}$ to be the smallest
\mbox{$\bigl( \tilde A - \tilde B_{K} (V_K \tilde D_K)^\dag V_K \tilde
  C \,,\, \Ker(Q_{K} \tilde C) \bigr)$-conditioned invariant} subspace
containing $\Image(\bar B_K)$, and let $J_K$ be an output injection
matrix such that
%(see \cite[Chapter 5]{HLT-AS-MH:01})
\begin{align}
  \label{eq: S star}
  (\tilde A - \tilde B_{K} (V_K \tilde D_K)^\dag V_K \tilde C +J_K
  Q_{K} \tilde C) \Star \subseteq \Star.
\end{align}
Let $P_{K}$ be an orthonormal projection matrix onto the quotient
space $\real^{2n} \setminus \Star$, and let 
\begin{align}
  \label{eq: Ak}
  A_{K} = P_{K} (\tilde A - \tilde B_{K} (V_K \tilde D_K)^\dag V_K
  \tilde C + J_{K} Q_{K} \tilde C)P_{K}^{\transpose}.
\end{align}
Finally, let $H_K$ and the unique $M_K$ be such that
\begin{align}
  \label{eq:definition H and M}
  \begin{split}
    \Ker(H_K Q \tilde C) &= \Star + \Ker(Q \tilde C),\text{ and }\\
    H_K Q \tilde C &= M_K P_{K}.
  \end{split}
\end{align}

% \begin{enumerate}

% 	\item For any $k$-dimensional subset $K$ of $\until{2n+m+p}$, let $[V_{K}^\transpose \; Q_{K}^\transpose]^\transpose \in \mathbb R^{p \times p}$ be
%   an orthonormal matrix, with $V_{K} = \Basis(\Image(\tilde
%   D_K))$ and $Q_{K} = \Basis(\Image(\tilde D_K)^\perp)$.

% 	\item Let $B_{Z} = \tilde B_{K} (V_{K} \tilde D_K)^\dag$ and $\bar B_{K} = \tilde B_{K}(I-D_{K}D_{K}^{\dag})$.
	
% 	\item Define $\Star$ to be the smallest 
%   $\bigl( (\tilde A - \tilde B_{K} (V_K \tilde D_K)^\dag V_K \tilde C) \,,\, \Ker(Q_{K} \tilde C) \bigr)$-conditioned 
%   invariant containing $\Image(\bar B_K)$, and
  	
% 	\item Let $J_K$ be an output injection matrix such that $(\tilde A+J_K
%   Q_{K} \tilde C) \Star \subseteq \Star$.
  
%   	\item Let $P_{K}$ be an orthonormal projection matrix onto the quotient
%   space $\real^{2n} \setminus \Star$.
  
%   	\item Let $A_{K} = P_{K} (\tilde A - \tilde B_{K} (V_K \tilde D_K)^\dag V_K \tilde C + J_{K} Q_{K} \tilde C)P_{K}^{\transpose}$.
	
% 	\item Let $H_K$ be such that $\Ker(H_K Q \tilde C) = \Star +
%   \ker(Q \tilde C)$, and let $M_K$ be the unique solution to $H_K Q \tilde C = M_K P_{K}$.

% \end{enumerate}
  
\begin{theorem}[\bf Attack identification filter]
\label{Proposition: Attack identification filter}
Consider the power network descriptor system \eqref{eq:
  cyber_physical_fault} and the associated Kron-reduced system
\eqref{eq: power network Kron-reduced model}. Assume that the attack
set $K$ is identifiable and that the network initial state is known.
% Let the attack set $K$ be identifiable. Let $k \in \mathbb{N}$ be
% the cardinality of the attack set, and let $K$ be any
% $k$-dimensional subset of $\until{2n+m+p}$.
Consider the
\emph{identification filter}
  \begin{align}
  \label{eq: identification filter}
  \begin{split}
    \dot w_{K} (t) 
    =& 
    (A_{K} + G_{K} M_{K}) w_{K} (t) 
    \\ &
    + \bigl( P_K B_{Z} V_{K} - (P_{K} J_{K} + G_{K} H_{K}) Q \bigr) y (t)
    ,\\
    r_{K} (t) 
    =& 
    M_{K} w_{K} (t) - H_{K} Q y(t),
    \end{split}
  \end{align}
  where $w_{K} (0) = P_{K} x(0)$, and $G_K \in \mathbb R^{2n \times
    p}$
  is such that $A_K + G_K M_K$ is a Hurwitz matrix. Then $r_{K} (t) =
  0$ at all times $t \in \real_{\ge 0}$ if and only if $K$ equals the
  attack set.
  \end{theorem}
  \smallskip

  Note that the residual $r_K(t)$ is identically zero if the attack
  set coincides with $K$, even if the attack input is nonzero.

  \begin{pfof}{Theorem \ref{Proposition: Attack identification filter}}
    Let $R$ be an attack set with $|R|\le|K|$ and $R \neq K$.
    % Let $K$ and $R$ be two $k$ dimensional attack sets of
    % $\until{2n+m+p}$, with $K \neq R$, and let $\tilde B_K$, $\tilde
    % D_K$, $\tilde B_{R}$, and $\tilde D_{R}$ be the associated input
    % matrices.
    % Let
    % $[V_{K}^\transpose \; Q_{K}^\transpose]^\transpose \in \mathbb
    % R^{p
    % \times p}$ be an orthonormal 
    With the output transformation
    % matrix , with
    % $V_{K} = \Basis(\Image(\tilde D_K))$ and $Q_{K} =
    % \Basis(\Image(\tilde D_K)^\perp)$. 
    % In the new output coordinates
    $[z_1,z_2] = [V_{K}y,Q_{K}y]$, the Kron-reduced system \eqref{eq:
      power network Kron-reduced model} becomes
  \begin{align}
    \label{eq: power network Kron-reduced model under attack - 1}
    \begin{split}
      \dot x (t) &= \tilde A x(t) + 
      \begin{bmatrix}
        \tilde B_K & \tilde B_{R}
      \end{bmatrix}
      \begin{bmatrix}
        u_K(t) \\ u_{R} (t)
      \end{bmatrix}
      ,\\
      z_1(t) &= V_K \tilde C x (t) + V_K \tilde D_K u_K(t) + V_K \tilde D_{R} u_{R}(t),\\
      z_2(t) &= Q_K \tilde C x (t) + Q_K \tilde D_{R} u_{R}(t).
    \end{split}
 \end{align}
 Note that the attack set $K$ affects only the output
 $z_{1}(t)$. The output equation for $z_{1}(t)$ can be solved for
 $u_{K}(t)$ as
 \begin{multline*}
   u_{K}(t)
   =
   (V_{K} \tilde D_{K})^{\dag} (z_{1}(t) - V_{K} \tilde C x(t)) 
   \\ - (V_{K} \tilde D_{K})^{\dag} V_{K} \tilde D_{R} u_{R}(t) + \subscr{u}{hom}(t)
   \,,
 \end{multline*}
 where $\subscr{u}{hom}(t) \in \Ker(V_{k} \tilde D_{K}) = \Ker( \tilde
 D_{K})$ and $u_{R}(t)$ are unknown signals, while $z_{1}(t)$ is
 known. The Kron-reduced system \eqref{eq: power network Kron-reduced
   model under attack - 1} can equivalently be written with unknown
 inputs $u_{K}(t)$ and $u_{R}(t)$, known input $z_{1}(t)$, and output
 $z_{2}(t)$ as
 \begin{align}
   \label{eq: power network Kron-reduced model under attack - 2}
   \begin{split}
     \dot x (t) =& 
    (\tilde A - \tilde B_{K} (V_K \tilde D_K)^\dag V_K \tilde C) x(t) 
    \\&+ 
    \begin{bmatrix}
      B_Z & \bar B_{K} & \bar B_{R}
    \end{bmatrix}
    \begin{bmatrix}
      z_1(t) \\ u_K(t) \\ u_{R} (t)
    \end{bmatrix}
    ,\\
    z_2(t) =& Q_K \tilde C x (t) + Q_K \tilde D_{R} u_{R}(t) ,
    \end{split}
 \end{align}
 where $B_Z$ and $\bar B_K$ are as in \eqref{eq: Bz Bk}, and
 \begin{align*}
   \bar B_{R} = \tilde B_{R} - \tilde B_{K} (V_{K} \tilde
   D_{K})^{\dag} V_{K} \tilde D_{R}.
 \end{align*}
 % $B_{Z} = \tilde B_{K} (V_{K} \tilde
 % D_K)^\dag$, $\bar B_{K} = \tilde B_{K}(I-D_{K}D_{K}^{\dag})$, and
 % $\bar B_{R} = \tilde B_{R} - \tilde B_{K} (V_{K} \tilde D_{K})^{\dag}
 % V_{K} \tilde D_{R}$.
 Let $\Star$ and $J_K$ be as in \eqref{eq: S star},
 % In order to continue, let $\Star$ be the smallest $\bigl( (\tilde A -
 % \tilde B_{K} (V_K \tilde D_K)^\dag V_K \tilde C) \,,\, \Ker(Q_{K}
 % \tilde C) \bigr)$-conditioned invariant containing $\Image(\bar B_K)$
 % \cite{HLT-AS-MH:01}. Recall that, for the system \eqref{eq: power
 %   network Kron-reduced model} there exists an observer for the
 % quotient space $\real^{2n+m} \setminus \Star$, which is nontrivial
 % due to attackers identifiability.  \flomargin{why non-trivial?}
 % Equivalently, there exists an output injection matrix $J_K$ such that
 % $\bigl( \tilde A - \tilde B_{K} (V_K \tilde D_K)^\dag V_K \tilde C +
 % J_K Q_{K} \tilde C \bigr) \Star \subseteq \Star$. 
 and consider the orthonormal change of coordinates given by $T_{K} =
 [W_K^\transpose \; P_K^\transpose] \in \mathbb R^{2n \times 2n}$,
 where $W_K^\transpose$ is a basis of $\Star$, $P_K$ is a projection
 matrix onto the quotient space $\real^{2n} \setminus \Star$, and
 $T_{K}^{-1} = T_{K}^\transpose$.  In the new coordinates
 $[\xi_1,\xi_2] = [W_{K}x,P_{K}x]$, system \eqref{eq: power network
   Kron-reduced model under attack - 2} reads as
 \begin{align}
 \label{eq: power network Kron-reduced model under attack - 3}
 \begin{split}
   \begin{bmatrix}
     \dot \xi_1 \\ \dot \xi_2
   \end{bmatrix}
   &\!=\!
   \begin{bmatrix}
     \hat A_{11} & \hat A_{12}\\
     0       & \hat A_{22}
   \end{bmatrix}\!
    \begin{bmatrix}
      \xi_1 \\ \xi_2
    \end{bmatrix}
    \!+\!\!
    \begin{bmatrix}
      \hat B_{11} & \hat B_{12} & \hat B_{13}\\
      \hat B_{21} & 0 & \hat B_{23}
    \end{bmatrix}\!
    \begin{bmatrix}
      z_1 \\u_K \\ u_{R}
    \end{bmatrix}
    \!,\\
    z_{2}(t)
    &\!=\!
    \begin{bmatrix}
      \hat C_{1} & \hat C_{2}
    \end{bmatrix}\!
    \begin{bmatrix}
      \xi_1 \\ \xi_2
    \end{bmatrix}
    \!+\! 
    \hat D u_{R}(t)
    \!.
    \end{split}
  \end{align}
  The zero pattern in the system and input matrix of \eqref{eq: power
    network Kron-reduced model under attack - 3} arises due to the
  invariance properties of $\Star$, which contains $\Image(\bar B_K)$.
  % An investigation of the matrix pencil of system \eqref{eq: power
  %   network Kron-reduced model under attack - 3} shows that, under the
  % identifiability hypothesis, the subsystem $(\hat A_{22},\hat
  % B_{23},\hat C_{2},\hat D)$ has only trivial zero dynamics.
  %
  % Likewise, since $\Star$ is contained in the subspace reachable
  % from $u_K$, the pair $(\hat A_{11}, \hat B_{11})$ is completely
  % reachable. Finally, notice that the subsystem $(\hat A_{22},\hat
  % B_{22},\hat C_{2},\hat D)$ has only trivial zero dynamics, because
  % of the identifiability hypothesis, and because output injection
  % does not alter the invariant zeros of a linear system.
  For the system \eqref{eq: power network Kron-reduced model
    under attack - 3} we propose the filter
  \begin{align}
    \dot w_{K} (t) &\!=\! (\hat A_{22} + G_{K} M_{K})
    w_{K} (t) + \hat B_{21} z_1(t) - G_{K} H_{K} z_2 (t),
    \nonumber\\
    r_{K} (t) &\!=\! M_{K} w_{K} (t) - H_{K} z_2(t),
    \label{eq: identification filter - nice coordinates}
  \end{align}
  where $G_K$ is chosen such that $\hat A_{22} + G_{K} M_{K}$ is a
  Hurwitz matrix. Let $H_K$ and $M_K$ be as in \eqref{eq:definition H
    and M}, which, in these coordinates, coincides with $H_K \hat C_1
  = 0$ and $H_K \hat C_2 = M_K$. Define the filter error $e(t) =
  w_K(t) - \xi_2(t)$, then the residual filter \eqref{eq:
    identification filter - nice coordinates} written in error
  coordinates is
  \begin{align*}
    \dot e(t) =& (\hat A_{22} + G_{K} M_{K}) w_{K} (t) -\hat A_{22} \xi_2 (t)
    - \hat B_{23} u_{R} (t)\\ &- G_{K} H_{K} \big([ \hat C_{1} \; \hat C_{2}]
    \xi(t) + \hat D u_{R}(t) \big) \\
    =& (\hat A_{22} + G_{K} M_{K}) e(t) - (\hat B_{23} + G_{K} H_{K} \hat D) u_{R}(t) \\
    r_K(t) =& M_K e(t) - H_K \hat D u_{R}(t).
  \end{align*}
  It can be shown that $(\hat A_{22} + G_{K} M_{K}, - (\hat B_{23} +
  G_{K} H_{K} \hat D), M_K, - H_K \hat D )$ has no zero dynamics, so
  that the residual $r_K (t)$ is not affected by $K$, and every
  nonzero signal $u_R (t)$ is detectable from $r_K (t)$. Consequently,
  $r_K(t)$ is identically zero if and only if $K$ is the attack
  set. Finally, in original coordinates, the filter \eqref{eq:
    identification filter - nice coordinates} takes the form
  \eqref{eq: identification filter}.
  % Because output injection operations do not alter the invariant zeros
  % of a system, we conclude that the error system $(\hat A_{22} + G_{K}
  % M_{K}, -\hat B_{23} - G_K H_K \hat D, M_K, -H_K \hat D)$ has only
  % trivial zero dynamics. Note also that $u_K(t)$ does not affect
  % $r_K(t)$. Consequently, $r_K(t)$ is identically zero if and only if
  % $K$ is the attack set.
\end{pfof}

For an attack set $K$, we refer to the signal $r_K (t)$ in the
filter~\eqref{eq: identification filter} as the residual associated
with $K$. A corollary result of Theorem \ref{Proposition: Attack
  identification filter} is that, if only an upper bound on the
cardinality of the attack set is known, then the residual $r_K (t)$ is
nonzero if and only if the attack set is contained in $K$. We now
summarize our identification procedure, which assumes the knowledge of
the network initial condition and of an upper bound $k$ on the
cardinality of the attack set $K$:
\begin{enumerate}
\item design an identification filter for each possible subset of
  $\until{2n+m+p}$ of cardinality $k$;
\item monitor the power network by running each identification filter;
\item the attack set $K$ coincides with the intersection of the attack
  sets $Z$ whose residual $r_Z (t)$ is identically zero.
\end{enumerate}

\begin{remark}\emph{(\bf Detection and identification filters for unknown initial condition)}
\label{Remark: detection filter for unknown initial condition}
% The detection and the identification filters in Theorem
% \ref{Proposition: Attack detection filter} and \ref{Proposition:
% Attack identification filter} require the knowledge of the network
% initial state.
% For the implementation of the detection and identification filters,
If the network initial state is not available, then an arbitrary
initial state $w(0) \in \mathbb R^{2n}$ can be chosen. Consequently,
the filters performance becomes asymptotic, and some attacks may
remain undetected or unidentified. For instance, if the eigenvalues of
the detection filter matrix have been assigned to have real part
smaller than $c < 0$, with $c \in \real$, then, in the absence of
attacks, the residual $r(t)$ exponentially converges to zero with rate
less than $c$. Hence, only inputs $u(t)$ that vanish faster or equal
than $e^{-ct}$ can remain undetected by the filter \eqref{eq:
  detection filter}. Alternatively, the detection filter can be
modified so as to converge in a predefined finite time
\cite{AVM-HTT:94}. In this case, every attack signal is detectable
after a finite transient.
%TR-FA:08
%If such a design is used, then 
% attack detection can be performed exactly after a predefined initial
% transient. 
\oprocend
% and the filter \eqref{eq: detection filter} can be modified in
% either of the following two ways. As first modification, let $G$ be
% such that all eigenvalues of $\tilde A + G \tilde C$ have real part
% less than some constant $c < 0$, then the filter output $r(t)$
% exponentially converges to zero with rate less than $c$ if and only
% if\fpmargin{Sure of this? Certainly necessary, but also sufficient?}
% $u(t)$ converges exponentially to zero with rate less than $c$.
% Hence, any input $u(t)$ that does not vanish faster or equal than
% $e^{-ct}$ can be detected by the filter \eqref{eq: detection
% filter}. Alternatively, the filter \eqref{eq: detection filter} can
% be implemented such that it is convergent in a finite time interval
% of length $\tau > 0$ by using the modifications suggested, for
% instance, in \cite{AVM-HTT:94,TR-FA:08}. Hence, for all $t > \tau$
% we have that $r(t) = 0$ if and only if $u(t) = 0$ and non-zero
% inputs can be detected.
\end{remark}

\begin{remark}\emph{(\bf Detection and identification in the presence
    of process and measurement noise)}
  The detection and identification filters here presented are a
  generalization to dynamical systems with direct input to output
  feedthrough of the devices presented in
  \cite{MAM-GCV-ASW-CM:89}. Additionally, our design guarantees the
  absence of invariant zeros in the residual system, so that
  \emph{every} attack signal affect the corresponding
  residual. Finally, if the network dynamics are affected by noise,
  then an optimal noise rejection in the residual system can be
  obtained by choosing the matrix $G$ in \eqref{eq: detection filter}
  and $G_K$ in \eqref{eq: identification filter} as the Kalman gain
  according to the noise statistics. \oprocend
\end{remark}

% \begin{remark}\emph{(\bf Identification filter for unknown initial condition)}
% \label{Remark: identification filter for unknown initial condition}
% The filter \eqref{eq: identification filter} in Proposition
% \ref{Proposition: Attack identification filter} requires the knowledge
% of the initial state $x(0)$. If the initial state $x(0)$ is not known,
% then, analogously to the detection filter \eqref{eq: detection
%   filter}, the identification filter can be modified as stated in
% Remark \ref{Remark: detection filter for unknown initial condition}.
% \oprocend
% \end{remark}

% Notice that for a successful identification of all possible attack
% sets of cardinality $k$, Proposition \ref{Proposition: Attack
%   identification filter} requires one filter for each attack. Hence,
% overall $\left(\begin{smallmatrix}
%   2n+m+p \\ 2
% \end{smallmatrix}\right)$
% filters are needed.

%%% Local Variables: 
%%% mode: latex
%%% TeX-master: "main"
%%% End: 

\section{A numerical study}\label{sec:example}
The effectiveness of our theoretic developments is here demonstrated
for the IEEE 14 bus system reported in Fig. \ref{fig:ieee14}.
Let the IEEE 14 bus power network be modeled as a descriptor model of
the form \eqref{eq: cyber_physical_fault}, where the network matrix
$A$ is as in \cite{rdz-cem-dg:11}. Following \cite{YL-MKR-PN:09}, the
measurement matrix $C$ consists of the real power injections at all
buses, of the real power flows of all branches, and of one rotor angle
(or one bus angle). We assume that an attacker can independently
compromise every measurement, except for the one referring to the
rotor angle, and that it does not inject state attacks.
% and only the network measurements. Hence, we let $F$ equal the zero
% matrix, and $L$ equal the identity matrix.

Let $k \in \mathbb{N}$ be the cardinality of the attack set. From
\cite{YL-MKR-PN:09} it is known that, for a Static Detector, an
undetectable attack exists if $k \ge 4$. In other words, due to the
sparsity pattern of $C$, there exists a signal $u_K(t)$, with (the
same) four nonzero entries at all times, such that $D u_K(t) \in
\Image (C)$ at all times. By Theorem \ref{Theorem: Static
  detectability of cyber-physical attacks} the attack set $K$ remains
undetected by a Static Detector through the attack mode $u_K (t)$. On
the other hand, following Theorem \ref{Theorem: Dynamic detectability
  of cyber-physical attacks}, it can be verified that, for the same
output matrix $C$, and independent of the value of $k$, there exists
\emph{no} undetectable (output) attack set.

\section{Conclusion}\label{sec:conclusion}
For a power network modeled via a linear time-invariant descriptor
system, we have analyzed the fundamental limitations of static and
dynamic attack detection and identification procedures. We have
rigorously shown that a dynamic detection and identification method
exploits the network dynamics and outperforms the static counterpart,
while requiring, possibly, fewer measurements. Additionally, we have
described a provably correct attack detection and identification
procedure based on dynamic residuals filters, and we have illustrated
its effectiveness through an example of cyber-physical attacks against
the IEEE 14 bus system.

%%% Local Variables: 
%%% mode: latex
%%% TeX-master: "main"
%%% End:

% \section{Proof no zero dynamics}
% \input{proof_no_zero_dynamics}

% \begin{appendix}

% \input{modelling}

% \end{appendix}

%%%%%%%%%%%%%%%%%%%%%%%%%%%%%%%%%%%%%%%%%%%%%%%%%%%%%%%%%%%%%%%%%%%%%%%%%%%%%%%%

%\renewcommand{\baselinestretch}{0.97}

\bibliographystyle{ieeeconf}

\end{document}